# Discrete chain graph models

MATHIAS DRTON

*Department of Statistics, University of Chicago, Chicago, IL 60637, USA.*
*E-mail: drton@uchicago.edu*

The statistical literature discusses different types of Markov properties for chain graphs that lead to four possible classes of chain graph Markov models. The different models are rather well understood when the observations are continuous and multivariate normal, and it is also known that one model class, referred to as models of LWF (Lauritzen–Wermuth–Frydenberg) or block concentration type, yields discrete models for categorical data that are smooth. This paper considers the structural properties of the discrete models based on the three alternative Markov properties. It is shown by example that two of the alternative Markov properties can lead to non-smooth models. The remaining model class, which can be viewed as a discrete version of multivariate regressions, is proven to comprise only smooth models. The proof employs a simple change of coordinates that also reveals that the model's likelihood function is unimodal if the chain components of the graph are complete sets.

## 1. Introduction

A graphical Markov model is a statistical model defined over a graph whose vertices correspond to observed random variables. The missing edges of the graph are translated into conditional independence restrictions that the model imposes on the joint distribution of the variables [21]. Among the more complex graphical models are those based on chain graphs. Chain graphs may have both directed and undirected edges under the constraint that there do not exist any semi-directed cycles. The absence of semi-directed cycles implies that the vertex set of a chain graph can be partitioned into so-called chain components such that edges within a chain component are undirected whereas the edges between two chain components are directed and point in the same direction.

The rules that govern how a graph is translated into conditional independence restrictions are known as Markov properties. Four classes of Markov properties for chain graphs have been discussed in the literature, and we classify them as:

Type I: the LWF or block concentration Markov property [12, 22];

Type II: the AMP (alternative Markov property) or concentration regression Markov property [1];
Type III: a Markov property that is dual to the type II property;
Type IV: the multivariate regression Markov property [4, 27], which can also be viewed as a special case of Markov properties for path diagrams [20, 25, 26].

The four types arise by combining two different interpretations of directed edges with two different interpretations of undirected edges (compare Section 2).

The four classes of Gaussian (i.e., multivariate normal) chain graph models associated with the above Markov properties are rather well understood. In particular, they are known to be smooth (i.e., they are curved exponential families [17]). Discrete models for categorical data have been thoroughly explored under the Markov property of type I (LWF). The resulting models have log-linear structure [21], Section 4.6.1, which yields that the models are smooth exponential families. However, less is known about the other discrete models. We note that discrete type IV models are related to models employed in longitudinal data analysis; see, for example, [11], and that, despite being termed modified path diagram models, the models discussed in [14] are most closely related to type I models.

This paper investigates smoothness properties of discrete models of type II, III and IV, which are all algebraic exponential families in the sense of [9]. Studying smoothness is important because the standard asymptotic distribution theory (e.g., normal distribution limits for maximum likelihood estimators and $\chi^2$-limits for likelihood ratios) is valid in smooth algebraic exponential families but may fail in non-smooth models [6]. Smoothness of conditional independence models cannot be taken for granted as demonstrated by the following example; see also [2], Example 7. The example concerns two discrete random variables $X_1$ and $X_2$ that are independent marginally as well as conditionally given a third binary variable $X_3$; in symbols, $X_1 \perp\!\!\!\perp X_2$ and $X_1 \perp\!\!\!\perp X_2 \mid X_3$. The corresponding subset of the appropriate probability simplex is a union of two sets corresponding to $X_1 \perp\!\!\!\perp (X_2, X_3)$ and $X_2 \perp\!\!\!\perp (X_1, X_3)$, respectively. The set defined by $X_1 \perp\!\!\!\perp (X_2, X_3)$ is a smooth manifold and so is the set given by $X_2 \perp\!\!\!\perp (X_1, X_3)$. Their union, however, fails to be smooth where the two components intersect. This intersection corresponds to complete independence of $X_1$, $X_2$ and $X_3$. Details on how the presence of singularities in this example affects the behavior of a likelihood ratio test can be found in [9], Section 4.2 and [6], Example 2.7, where the Gaussian version of the problem is treated. The Gaussian case is analogous since for a jointly multivariate normal random vector $(X_1, X_2, X_3)$ it also holds that $X_1 \perp\!\!\!\perp X_2$ and $X_1 \perp\!\!\!\perp X_2 \mid X_3$ is equivalent to $X_1 \perp\!\!\!\perp (X_2, X_3)$ or $X_2 \perp\!\!\!\perp (X_1, X_3)$.

The main result of this paper shows that discrete type IV models are smooth. Stated in Corollary 10, this result follows from a linear change of conditional probability coordinates that simplifies the conditional independence constraints in the model definition (Theorem 8 in Section 3). Moreover, type IV models have unimodal likelihood functions if the chain components of the underlying graph are complete sets, in which case the models of type II and type IV coincide (Section 4). Finally, we show by example in Sections 5 and 6 that the classes of type II and III include non-smooth models. The paper concludes with the discussion in Section 7.

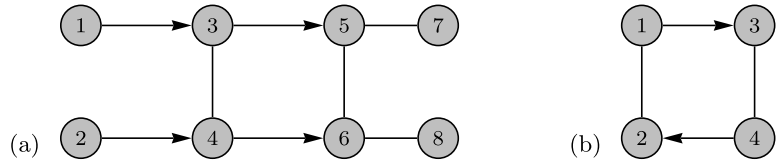

**Figure 1.** (a) Chain graph with chain components $\{1\}$, $\{2\}$, $\{3,4\}$ and $\{5,6,7,8\}$; (b) a graph that is not a chain graph.

## 2. Chain graphs and their Markov properties

### 2.1. Chain graphs

Let $G = (V, E)$ be a graph with finite vertex set $V$ and edge set $E \subseteq (V \times V) \setminus \{(v,v) \mid v \in V\}$. An edge $(v, w) \in E$ is *directed* if $(w, v) \notin E$ and *undirected* if $(w, v) \in E$. We denote a directed edge $(v, w)$ by $v \to w$ and write $v - w$ if $(v, w)$ is undirected. If $(v, w) \in E$ then $v$ and $w$ are *adjacent*. If $v \to w$ then $v$ is a *parent* of $w$, and if $v - w$ then $v$ is a *neighbor* of $w$. Let $\mathrm{pa}_G(v)$ and $\mathrm{nb}_G(v)$ denote the sets of parents and neighbors of $v$, respectively. For $\sigma \subseteq V$, let

$$\mathrm{pa}_G(\sigma) = \left(\bigcup_{v \in \sigma} \mathrm{pa}_G(v)\right) \setminus \sigma,$$

$$\mathrm{nb}_G(\sigma) = \left(\bigcup_{v \in \sigma} \mathrm{pa}_G(v)\right) \setminus \sigma,$$

and $\mathrm{Nb}_G(\sigma) = \mathrm{nb}_G(\sigma) \cup \sigma$.

A *path* in $G$ is a sequence of distinct vertices $\langle v_0, \ldots, v_k \rangle$ such that $v_{i-1}$ and $v_i$ are adjacent for all $1 \leq i \leq k$. A path $\langle v_0, \ldots, v_k \rangle$ is a *semi-directed cycle* if $(v_i, v_{i+1}) \in E$ for all $0 \leq i \leq k$ and at least one of the edges is directed as $v_i \to v_{i+1}$. Here, $v_{k+1} \equiv v_0$. A *chain graph* is a graph without semi-directed cycles (see Figure 1). Define two vertices $v_0$ and $v_k$ in a chain graph $G$ to be equivalent if there exists a path $\langle v_0, \ldots, v_k \rangle$ such that $v_i - v_{i+1}$ in $G$ for all $0 \leq i \leq k-1$. The equivalence classes under this equivalence relation are the *chain components* of $G$. The chain components $(\tau \mid \tau \in \mathcal{T})$ yield a partitioning of the vertex set

$$V = \dot{\bigcup_{\tau \in \mathcal{T}}} \tau, \tag{2.1}$$

and the subgraph $G_\tau$ induced by each chain component $\tau$ is a connected undirected graph. Moreover, the directed edges between two chain components $\tau_1$ and $\tau_2$ all have the same direction, that is, if $(v, w) \in \tau_1 \times \tau_2$ and $(x, y) \in \tau_1 \times \tau_2$ are two pairs of adjacent vertices, then either $v \to w$ and $x \to y$ or $w \to v$ and $y \to x$ in $G$. It follows that we can define an acyclic digraph (DAG) $D = D(G)$ over the chain components: $\mathcal{T}$ is the vertex set of $D$, and we draw an edge $\tau_1 \to \tau_2$ if and only if there exists $v \in \tau_1$ and $w \in \tau_2$ with $v \to w$ in $G$.

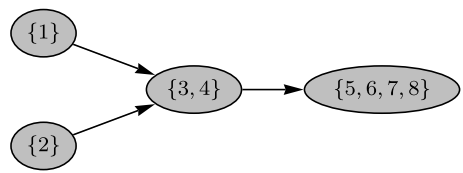

**Figure 2.** DAG of chain components for the chain graph from Figure 1(a).

**Example 1.** Consider the chain graph $G$ in Figure 1(a). It has four chain components, namely, $\{1\}, \{2\}, \{3,4\}$ and $\{5,6,7,8\}$. The DAG $D(G)$ has these four chain components as nodes and is depicted in Figure 2.

### 2.2. Block-recursive Markov properties

A *Markov property* of a graph $G = (V, E)$ lists conditional independence statements

$$\alpha_i \perp\!\!\!\perp \beta_i \mid \gamma_i, \qquad i = 1, \ldots, k,$$

for triples $(\alpha_i, \beta_i, \gamma_i)$ of pairwise disjoint subsets of $V$ with $\alpha_i, \beta_i \neq \varnothing$. These triples are determined by the edge set $E$. The joint distribution $P$ of a random vector $X \in \mathbb{R}^V$ *obeys* the Markov property if for all $1 \leq i \leq k$, the subvector $X_{\alpha_i}$ is conditionally independent of $X_{\beta_i}$ given $X_{\gamma_i}$. If $\gamma_i = \varnothing$ then the conditional independence is understood as marginal independence of $X_{\alpha_i}$ and $X_{\beta_i}$.

Block-recursive Markov properties for a chain graph $G$ employ the recursive structure of the chain components captured in the DAG $D = D(G)$; see [1, 21]. For $\tau \in \mathcal{T}$, let $\mathrm{pa}_D(\tau)$ be the union of all $\bar{\tau} \in \mathcal{T} \setminus \{\tau\}$ that are parents of $\tau$ in $D$. Similarly, the set of *non-descendants* $\mathrm{nd}_D(\tau)$ is the union of all $\bar{\tau} \in \mathcal{T} \setminus \{\tau\}$ for which there is no directed path $\bar{\tau} \to \cdots \to \tau$ in $D$. The following conditional independence statements are associated with the DAG $D$:

$$\tau \perp\!\!\!\perp [\mathrm{nd}_D(\tau) \setminus \mathrm{pa}_D(\tau)] \mid \mathrm{pa}_D(\tau) \qquad \forall \tau \in \mathcal{T}. \tag{C1}$$

If the joint distribution has a density with respect to a product measure, then the conditional independence relations in (C1) are equivalent to the density factorizing over the graph [21]. We will employ the factorization over the DAG $D$ in our study of discrete chain graph models; see (3.2).

**Example 1 (Cont.).** If we consider the chain component $\tau = \{1\}$ of the chain graph $G$ from Figure 1(a), then $\mathrm{pa}_D(\tau) = \varnothing$ and $\mathrm{nd}_D(\tau) = \{2\}$. Hence, (C1) states

$$\{1\} \perp\!\!\!\perp \{2\}.$$

If $\tau = \{5,6,7,8\}$, then $\mathrm{pa}_D(\tau) = \{3,4\}$ and $\mathrm{nd}_D(\tau) = \{1\} \cup \{2\} \cup \{3,4\}$, which leads to

$$\{5,6,7,8\} \perp\!\!\!\perp \{1,2\} \mid \{3,4\}.$$

The factorization of a joint density $f(x_1, \ldots, x_7)$ alluded to above takes the form

$$f(x_1, \ldots, x_8) = f(x_1) f(x_2) f(x_3, x_4 \mid x_1, x_2) f(x_5, x_6, x_7, x_8 \mid x_3, x_4).$$

Each chain component $\tau \in \mathscr{T}$ induces an undirected subgraph $G_\tau$. Applying the local version of the classic Markov property for these undirected graphs (see, e.g., [21]) to the conditional distribution of the variables in $\tau$ given the variables in $\mathrm{pa}_D(\tau)$ leads to the conditional independence statements

$$\sigma \perp\!\!\!\perp [\tau \setminus \mathrm{Nb}_G(\sigma)] \mid [\mathrm{pa}_D(\tau) \cup \mathrm{nb}_G(\sigma)] \qquad \forall \tau \in \mathscr{T}, \forall \sigma \subseteq \tau. \tag{C2a}$$

If the conditional distributions for $\tau$ given $\mathrm{pa}_D(\tau)$ have positive densities with respect to a product measure, then the Hammersley–Clifford theorem implies that the conditional independence relations in (C2a) correspond to factorizations of the conditional densities; see again [21] for the precise results. As an alternative to (C2a), we can employ a dual Markov property for undirected graphs (discussed, e.g., in [18]) that yields

$$\sigma \perp\!\!\!\perp [\tau \setminus \mathrm{Nb}_G(\sigma)] \mid \mathrm{pa}_D(\tau) \qquad \forall \tau \in \mathscr{T}, \forall \sigma \subseteq \tau. \tag{C2b}$$

Both (C2a) and (C2b) describe the consequences of the absence of undirected edges within a chain component $\tau$, but contrary to (C2a), the conditional independence relations in (C2b) are generally not related to density factorizations.

*Example 1 (Cont.).* Let $\tau = \{5,6,7,8\}$ be the largest chain component of the graph $G$ in Figure 1(a), for which $\mathrm{pa}_D(\tau) = \{3,4\}$. If $\sigma = \{5,7\}$, then $\mathrm{nb}_G(\sigma) = \{6\}$ and $\mathrm{Nb}_G(\sigma) = \{5,6,7\}$. Therefore, (C2a) states that

$$\{5,7\} \perp\!\!\!\perp \{8\} \mid \{3,4,6\},$$

whereas (C2b) states that

$$\{5,7\} \perp\!\!\!\perp \{8\} \mid \{3,4\}.$$

The final ingredient to the block-recursive Markov properties describes finer dependence structures associated with the absence of directed edges. Again there are two versions, namely,

$$\sigma \perp\!\!\!\perp [\mathrm{pa}_D(\tau) \setminus \mathrm{pa}_G(\sigma)] \mid [\mathrm{pa}_G(\sigma) \cup \mathrm{nb}_G(\sigma)] \qquad \forall \tau \in \mathscr{T}, \forall \sigma \subseteq \tau \tag{C3a}$$

and

$$\sigma \perp\!\!\!\perp [\mathrm{pa}_D(\tau) \setminus \mathrm{pa}_G(\sigma)] \mid \mathrm{pa}_G(\sigma) \qquad \forall \tau \in \mathscr{T}, \forall \sigma \subseteq \tau. \tag{C3b}$$

The two versions differ by whether vertices from the considered chain component $\tau$ are included in the conditioning set or not.

*Discrete chain graph models* 741

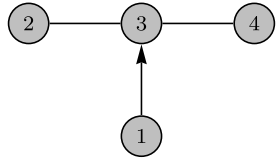

**Figure 3.** Chain graph with chain components $\{1\}$ and $\{2,3,4\}$.

*Example 1 (Cont.).* Consider again the graph $G$ in Figure 1(a) and the chain component $\tau = \{5,6,7,8\}$ with $\mathrm{pa}_D(\tau) = \{3,4\}$. If $\sigma = \{5,7\}$, then $\mathrm{pa}_G(\sigma) = \{3\}$ and $\mathrm{nb}_G(\sigma) = \{6\}$. Therefore, (C3a) states that

$$\{5,7\} \perp\!\!\!\perp \{4\} \mid \{3,6\},$$

whereas (C3b) states that

$$\{5,7\} \perp\!\!\!\perp \{4\} \mid \{3\}.$$

We are now ready to formally define the four types of Markov properties mentioned in the Introduction.

*Definition 2.* Let $G$ be a chain graph with chain components $(\tau \mid \tau \in \mathscr{T})$. A *block-recursive Markov property* for $G$ states (C1), one choice of either (C2a) or (C2b), and one choice of either (C3a) or (C3b). The block-recursive Markov property is of

$$type \left\{\begin{array}{c} I \\ II \\ III \\ IV \end{array}\right\} \text{ if it states (C1)}, \left\{\begin{array}{c} (\mathrm{C2a}) \\ (\mathrm{C2a}) \\ (\mathrm{C2b}) \\ (\mathrm{C2b}) \end{array}\right\} \text{ and } \left\{\begin{array}{c} (\mathrm{C3a}) \\ (\mathrm{C3b}) \\ (\mathrm{C3a}) \\ (\mathrm{C3b}) \end{array}\right\}.$$

As can be seen also in the next example, the Markov property of type I is the 'most conditional' with the largest conditioning sets, whereas type IV is the 'most marginal' with the smallest conditioning sets. Types II and III mix marginal and conditional perspectives.

*Example 3.* Let $G$ be the chain graph in Figure 3. Then (C1) is a void statement. The remaining statements can be summarized as follows:

$$\begin{aligned} \text{type I:} \quad & 2 \perp\!\!\!\perp 4 \mid \{1,3\} \text{ and } 1 \perp\!\!\!\perp \{2,4\} \mid 3, \\ \text{type II:} \quad & 2 \perp\!\!\!\perp 4 \mid \{1,3\} \text{ and } 1 \perp\!\!\!\perp \{2,4\}, \\ \text{type III:} \quad & 2 \perp\!\!\!\perp 4 \mid 1 \text{ and } 1 \perp\!\!\!\perp \{2,4\} \mid 3, \\ \text{type IV:} \quad & 2 \perp\!\!\!\perp 4 \mid 1 \text{ and } 1 \perp\!\!\!\perp \{2,4\}. \end{aligned}$$

In Section 3 we will consider discrete chain graph models of type IV. When studying these models we will exploit the following useful simplification of (C2b) that is due to [25], Theorem 4. This simplification is based on *connected sets*, which are subsets of the vertex set that induce a connected subgraph.

**Lemma 4.** *A probability distribution obeys the conditional independence relations* (C2b) *if and only if it obeys*

$$\sigma \perp\!\!\!\perp [\tau \setminus \mathrm{Nb}_G(\sigma)] \mid \mathrm{pa}_D(\tau) \qquad \forall \tau \in \mathscr{T}, \forall \sigma \subseteq \tau, \sigma \text{ connected.} \qquad \text{(C2b-conn)}$$

***Remark 5.*** In a multivariate normal distribution, two pairwise marginal independences $v \perp\!\!\!\perp w$ and $v \perp\!\!\!\perp u$ imply that $v \perp\!\!\!\perp \{u,w\}$. As a consequence, Gaussian chain graph models can be discussed in terms of *pairwise Markov properties* that list conditional independences between pairs of non-adjacent vertices in the graph; compare, for example, [7, 27]. When considering discrete random vectors taking only finitely many values, the pairwise Markov property for type I models is still equivalent to the respective block-recursive property as long as one limits oneself to positive joint distributions [12], Theorem 3.3. However, for the models of types II/III/IV a positive discrete distribution that obeys the pairwise Markov property will generally not obey the block-recursive Markov property. This follows from the fact that in almost every positive joint distribution that exhibits $v \perp\!\!\!\perp w$ and $v \perp\!\!\!\perp u$, $v$ is not independent of both $\{u,w\}$. For instance, if $G$ is the graph in Figure 3, then the pairwise model of type II (AMP) would be based on $2 \perp\!\!\!\perp 4 \mid \{1,3\}$, $1 \perp\!\!\!\perp 2$ and $1 \perp\!\!\!\perp 4$. Similarly, the pairwise model of type IV (multivariate regression) would be based on $2 \perp\!\!\!\perp 4 \mid 1$, $1 \perp\!\!\!\perp 2$ and $1 \perp\!\!\!\perp 4$. Such pairwise models will generally be of considerably larger dimension than their block-recursive analogs.

## 3. Discrete models of type IV

Let $X = (X_v \mid v \in V)$ be a discrete random vector with component $X_v$ taking values in $[d_v] = \{1, \ldots, d_v\}$. Let $\mathcal{I} = \bigtimes_{v \in V}[d_v]$. For $i = (i_v \mid v \in V) \in \mathcal{I}$, let

$$p(i) = P(X_v = i_v \text{ for all } v \in V). \qquad (3.1)$$

The joint distribution of $X$ is determined by the probability vector

$$p = (p(i) \mid i \in \mathcal{I})$$

in the $|\mathcal{I}| - 1 = (\prod_{v \in V} d_v) - 1$ dimensional probability simplex $\Delta$. Let $\Delta^o$ be the interior of the probability simplex, which corresponds to the regular exponential family of positive distributions on $\mathcal{I}$.

***Definition 6.*** *The discrete chain graph model* $\mathbf{P}_{\mathrm{IV}}(G)$ *associated with the chain graph* $G = (V, E)$ *is the set of (positive) probability vectors in* $\Delta^o$ *that yield a distribution on* $\mathcal{I}$ *that obeys the block-recursive Markov property of type IV (multivariate regression).*

The higher-level structure of the model $\mathbf{P}_{\mathrm{IV}}(G)$ is determined by condition (C1). For a probability vector $p \in \Delta^o$, this condition is equivalent to the condition that

$$p(i) = \prod_{\tau \in \mathscr{T}} p(i_\tau \mid i_{\mathrm{pa}_D(\tau)}), \qquad i \in \mathcal{I}, \tag{3.2}$$

where

$$p(i_\tau \mid i_{\mathrm{pa}_D(\tau)}) = P(X_\tau = i_\tau \mid X_{\mathrm{pa}_D(\tau)} = i_{\mathrm{pa}_D(\tau)}).$$

The factorization in (3.2) is the usual factorization over a DAG, but applied to the DAG of chain components $D = D(G)$. For a subset $\alpha \subseteq V$, define $\mathcal{I}_\alpha = \bigtimes_{v \in \alpha} [d_v]$. For fixed $i_{\mathrm{pa}_D(\tau)} \in \mathcal{I}_{\mathrm{pa}_D(\tau)}$, the vector with components $p(i_\tau \mid i_{\mathrm{pa}_D(\tau)})$, $i_\tau \in \mathcal{I}_\tau$, is a probability vector in the interior of the $|\mathcal{I}_\tau| - 1$ dimensional probability simplex $\Delta_\tau^o$.

In [8], a linear change of coordinates is used to simplify the description of conditional independence constraints that correspond to (C2b). We will now show how to generalize this change of coordinates to a version involving conditional probabilities that simplifies both (C2b) and (C3b).

Consider a subset $\varnothing \neq \sigma \subseteq \tau$ of a chain component $\tau \in \mathscr{T}$. Let $i_{\mathrm{pa}_D(\tau)} \in \mathcal{I}_{\mathrm{pa}_D(\tau)}$ be a given conditioning state. Define the restricted state space

$$\mathcal{J}_\sigma = \bigtimes_{v \in \sigma} [d_v - 1].$$

The set $[d_v - 1] = \{1, \ldots, d_v - 1\}$ is the state space for random variable $X_v$ but with the highest-numbered state $d_v$ removed. (Any other state could be removed instead.) With each state $j_\sigma \in \mathcal{J}_\sigma$ we associate the conditional probability

$$q(j_\sigma | i_{\mathrm{pa}_D(\tau)}) = P(X_\sigma = j_\sigma | X_{\mathrm{pa}_D(\tau)} = i_{\mathrm{pa}_D(\tau)}), \tag{3.3}$$

which we call *saturated Möbius parameter* for $\sigma$ and $j_\sigma$ given $(\tau, i_{\mathrm{pa}_D(\tau)})$. For fixed chain component $\tau$ and conditional state $i_{\mathrm{pa}_D(\tau)}$, but $\sigma$ varying through the power set of $\tau$ and $i_\sigma$ varying through $\mathcal{J}_\sigma$, there are $|\mathcal{I}_\tau| - 1$ many saturated Möbius parameters. They can be computed from the conditional probabilities $p(i_\tau | i_{\mathrm{pa}_D(\tau)})$ by the obvious summations. These summations define a linear map $\mu : \Delta_\tau \to \mathbb{R}^{|\mathcal{I}_\tau| - 1}$ taking the conditional probabilities $p(i_\tau | i_{\mathrm{pa}_D(\tau)})$, $i_\tau \in \mathcal{I}_\tau$, to the saturated Möbius parameters given $(\tau, i_{\mathrm{pa}_D(\tau)})$. We note that these parameters are closely related to conditional versions of dependence ratios; see [10] and references therein.

The following fact corresponds to Proposition 6 in [8], which is based on the well-known Möbius inversion; see also [16] where the Kronecker product structure of the matrices for the linear maps $\mu$ and $\mu^{-1}$ is described.

**Lemma 7.** *The linear map $\mu : \Delta_\tau \to \mathbb{R}^{|\mathcal{I}_\tau| - 1}$ from conditional probabilities to saturated Möbius parameters is bijective with the inverse map determined as follows: Let $i_\tau \in \mathcal{I}_\tau$ and define*

$$\sigma = \sigma(i_\tau) := \{v \in \tau | i_v \in [d_v - 1]\} \subseteq \tau.$$

Let $q(j_\varnothing | i_{\mathrm{pa}_D(\tau)}) = 1$. Then

$$p(i_\tau | i_{\mathrm{pa}_D(\tau)}) = \sum_{\alpha : \sigma \subseteq \alpha \subseteq \tau} (-1)^{|\alpha \setminus \sigma|} \sum_{j_{\alpha \setminus \sigma} \in \mathcal{J}_{\alpha \setminus \sigma}} q(j_{\alpha \setminus \sigma}, i_\sigma | i_{\mathrm{pa}_D(\tau)}).$$

As we show next, both (C2b) and (C3b) take on a simple form when expressed in terms of the saturated Möbius parameter coordinates. We use that every set $\delta \subseteq \tau$ that is not connected in $G$ can be partitioned uniquely into inclusion-maximal connected sets $\gamma_1, \ldots, \gamma_r \subseteq \tau$,

$$\delta = \gamma_1 \dot\cup \gamma_2 \dot\cup \cdots \dot\cup \gamma_r. \tag{3.4}$$

**Theorem 8.** *Let $G$ be a chain graph with chain components $(\tau | \tau \in \mathcal{T})$. A probability vector $p \in \Delta^o$ belongs to the discrete chain graph model $\mathbf{P}_{\mathrm{IV}}(G)$ if and only if the following three conditions hold:*

(i) *The components of $p$ factor as in (3.2).*
(ii) *For all $\tau \in \mathcal{T}$ and $i_{\mathrm{pa}_D(\tau)} \in \mathcal{I}_{\mathrm{pa}_D(\tau)}$, the saturated Möbius parameters for $p$ satisfy that*

$$q(j_\delta | i_{\mathrm{pa}_D(\tau)}) = q(j_{\gamma_1} | i_{\mathrm{pa}_D(\tau)}) q(j_{\gamma_2} | i_{\mathrm{pa}_D(\tau)}) \cdots q(j_{\gamma_r} | i_{\mathrm{pa}_D(\tau)})$$

*for every disconnected set $\delta \subseteq \tau$ and $j_\delta \in \mathcal{J}_\delta$. Here $\gamma_1, \ldots, \gamma_r \subseteq \tau$ are the inclusion-maximal connected sets in (3.4).*

(iii) *For all $\tau \in \mathcal{T}$, connected subsets $\gamma \subseteq \tau$ and $j_\gamma \in \mathcal{J}_\gamma$, the saturated Möbius parameters for $p$ satisfy that*

$$q(j_\gamma | i_{\mathrm{pa}_D(\tau)}) = q(j_\gamma | k_{\mathrm{pa}_D(\tau)})$$

*for every pair $i_{\mathrm{pa}_D(\tau)}, k_{\mathrm{pa}_D(\tau)} \in \mathcal{I}_{\mathrm{pa}_D(\tau)}$ such that $i_{\mathrm{pa}_G(\gamma)} = k_{\mathrm{pa}_G(\gamma)}$.*

**Proof.** Clearly, the factorization (3.2) required in (i) is equivalent to (C1). Applying Theorem 8 in [8] to each of the conditional distributions associated with the different vectors $i_{\mathrm{pa}_D(\tau)} \in \mathcal{I}_{\mathrm{pa}_D(\tau)}$, we see that (ii) is equivalent to (C2b).

Condition (C3b) states that $\gamma \perp\!\!\!\perp \mathrm{pa}_D(\tau) \setminus \mathrm{pa}_G(\gamma) | \mathrm{pa}_G(\gamma)$ for all subsets $\gamma \subseteq \tau \in \mathcal{T}$ of the chain components $\mathcal{T}$. The conditional independence for given $\gamma$ and $\tau$ holds if and only if

$$p(i_\gamma | i_{\mathrm{pa}_D(\tau)}) = p(i_\gamma | k_{\mathrm{pa}_D(\tau)}) \tag{3.5}$$

for all $i_\gamma \in \mathcal{I}_\gamma$ and every pair $i_{\mathrm{pa}_D(\tau)}, k_{\mathrm{pa}_D(\tau)} \in \mathcal{I}_{\mathrm{pa}_D(\tau)}$ such that $i_{\mathrm{pa}_G(\gamma)} = k_{\mathrm{pa}_G(\gamma)}$. If $i_\gamma \in \mathcal{J}_\gamma$, then

$$p(i_\gamma | i_{\mathrm{pa}_D(\tau)}) = q(i_\gamma | i_{\mathrm{pa}_D(\tau)}). \tag{3.6}$$

Hence, (C3b) implies (iii).

For the reverse implication we claim that (ii) and (iii) imply (C3b). Fix $i_{\mathrm{pa}_D(\tau)}, k_{\mathrm{pa}_D(\tau)} \in \mathcal{I}_{\mathrm{pa}_D(\tau)}$ such that $i_{\mathrm{pa}_G(\gamma)} = k_{\mathrm{pa}_G(\gamma)}$. By Lemma 7 and (ii), we can express $p(i_\gamma | i_{\mathrm{pa}_D(\tau)})$ as a function of saturated Möbius parameters $q(j_\alpha | i_{\mathrm{pa}_D(\tau)})$, where $\alpha$ is a connected subset of $\gamma$ and $j_\alpha \in \mathcal{J}_\alpha$. By (iii), $q(j_\alpha | i_{\mathrm{pa}_D(\tau)}) = q(j_\alpha | k_{\mathrm{pa}_D(\tau)})$, and thus (3.5) holds. □

***Example 9.*** In the chain graph in Figure 3, chain component $\{1\}$ is a singleton without parents and the associated saturated Möbius parameters are $q(j_1)$, $j_1 \in [d_1 - 1]$. The saturated Möbius parameters for the second chain component $\{2, 3, 4\}$ are the conditional probabilities

$$q(j_2|i_1), \quad q(j_3|i_1), \quad q(j_4|i_1),$$
$$q(j_2, j_3|i_1), \quad q(j_2, j_4|i_1), \quad q(j_2, j_4|i_1), \quad q(j_2, j_3, j_4|i_1),$$

where $i_1 \in [d_1]$, $j_2 \in [d_2 - 1]$, $j_3 \in [d_3 - 1]$ and $j_4 \in [d_4 - 1]$. The saturated Möbius parameters correspond to a probability vector in $\mathbf{P}_{\mathrm{IV}}(G)$ if and only if the following equations hold for all $i_1, k_1 \in [d_1]$ and $j_v \in [d_v - 1]$, $v \geq 2$:

$$q(j_2, j_4|i_1) = q(j_2|i_1)q(j_4|i_1),$$
$$q(j_2|i_1) = q(j_2|k_1),$$
$$q(j_4|i_1) = q(j_4|k_1),$$
$$q(j_2, j_4|i_1) = q(j_2, j_4|k_1).$$

The first equation is given by Theorem 8(ii), the others by Theorem 8(iii).

Theorem 8 can be read as expressing certain saturated Möbius parameters as functions of the remaining ones. One obtains a parametrization of $\mathbf{P}_{\mathrm{IV}}(G)$ in which the parameters are the conditional probabilities

$$q_\gamma(j_\gamma | i_{\mathrm{pa}_G(\gamma)}) = P(X_\gamma = j_\gamma | X_{\mathrm{pa}_G(\gamma)} = i_{\mathrm{pa}_G(\gamma)}) \tag{3.7}$$

with $j_\gamma \in \mathcal{J}_\gamma$ and $i_{\mathrm{pa}_G(\gamma)} \in \mathcal{I}_{\mathrm{pa}_G(\gamma)}$. We call the probabilities in (3.7) the *Möbius parameters* for model $\mathbf{P}_{\mathrm{IV}}(G)$. For a chain component $\tau \in \mathscr{T}$, let $\mathcal{C}(\tau)$ be the family of connected sets in the induced subgraph $G_\tau$, and define the vector

$$q_\tau = (q_\gamma(j_\gamma | i_{\mathrm{pa}_G(\gamma)}) | \gamma \in \mathcal{C}(\tau), j_\gamma \in \mathcal{J}_\gamma, i_{\mathrm{pa}_G(\gamma)} \in \mathcal{I}_{\mathrm{pa}_G(\gamma)}).$$

Let $Q_\tau$ be the set of vectors $q_\tau$ that are obtained from some $p \in \mathbf{P}_{\mathrm{IV}}(G)$. Moreover, let $q = (q_\tau | \tau \in \mathscr{T})$ and define $Q_G$ to be the set of vectors $q$ obtained from some $p \in \mathbf{P}_{\mathrm{IV}}(G)$. The set

$$Q_G = \bigtimes_{\tau \in \mathscr{T}} Q_\tau \tag{3.8}$$

is a Cartesian product, that is, the Möbius parameters from different chain components are variation-independent. Each set $Q_\tau$, however, is constrained via polynomial inequalities and no additional factorization of $Q_\tau$ into a Cartesian product seems possible in general.

The equations among saturated Möbius parameters that appear in Theorem 8 are of particularly simple nature and reveal that the set $\mathbf{P}_{\mathrm{IV}}(G)$ is a smooth manifold in the interior of the probability simplex; compare, for example, [13], Theorem 1.

**Corollary 10.** *For every chain graph $G$ the discrete chain graph model $\mathbf{P}_{\mathrm{IV}}(G)$ is a curved exponential family.*

Corollary 10 implies that in discrete chain graph models of type IV, maximum likelihood estimators are asymptotically normal, likelihood ratio statistics for model comparisons have asymptotic $\chi^2$-distributions, and the Bayesian information criterion is consistent for model selection [17]. For application of likelihood ratio tests and information criteria, it is important to know the dimension of $\mathbf{P}_{\mathrm{IV}}(G)$, which is readily obtained from the chain graph $G$ and the numbers of states the involved random variables may take.

**Corollary 11.** *The dimension of the discrete chain graph model $\mathbf{P}_{\mathrm{IV}}(G)$ is*

$$\dim(\mathbf{P}_{\mathrm{IV}}(G)) = \sum_{\tau \in \mathscr{T}} \sum_{C \in \mathcal{C}(\tau)} \left( \prod_{v \in C} (d_v - 1) \right) \left( \prod_{w \in \mathrm{pa}_G(C)} d_w \right).$$

## 4. Likelihood inference in models of type IV

Continuing our discussion of models based on the multivariate regression Markov property (type IV), suppose we observe a sample $X^{(1)}, \ldots, X^{(n)}$ of independent and identically distributed discrete random vectors taking values in $\mathcal{I} = \bigtimes_{v \in V} [d_v]$. Suppose further that the joint distribution common to the random vectors in the sample corresponds to an unknown probability vector $p \in \mathbf{P}_{\mathrm{IV}}(G)$, where $G = (V, E)$ is a chain graph. If we define the counts

$$n(i) = \sum_{k=1}^{n} 1_{\{X^{(k)} = i\}}, \qquad i \in \mathcal{I},$$

then the *likelihood function* of $\mathbf{P}_{\mathrm{IV}}(G)$ is equal to

$$L(p) = \prod_{i \in \mathcal{I}} p(i)^{n(i)}. \tag{4.1}$$

This likelihood function admits a factorization that we express in (4.2) using the *log-likelihood function* $\ell(p) = \log L(p)$.

For $\alpha \subseteq V$ and $i_\alpha \in \mathcal{I}_\alpha$, let

$$n(i_\alpha) = \sum_{j \in \mathcal{I} : j_\alpha = i_\alpha} n(j)$$

and
$$p(i_\alpha) = P(X_\alpha = i_\alpha) = \sum_{j \in \mathcal{I}: j_\alpha = i_\alpha} p(j).$$

By (3.2), the log-likelihood function can be written as

$$\ell(p) = \sum_{\tau \in \mathcal{T}} \sum_{i_\tau \in \mathcal{I}_\tau} \sum_{i_{\text{pa}_D(\tau)} \in \mathcal{I}_{\text{pa}_D(\tau)}} n(i_\tau, i_{\text{pa}_D(\tau)}) \log p(i_\tau | i_{\text{pa}_D(\tau)}). \tag{4.2}$$

Since the Möbius parameter space $Q_G$ factors accordingly into a Cartesian product, see (3.8), it follows that we can maximize $\ell(p)$ over $\mathbf{P}_{\text{IV}}(G)$ by maximizing the *component log-likelihood functions*

$$\ell_\tau(p) = \sum_{i_\tau} \sum_{i_{\text{pa}_D(\tau)}} n(i_\tau, i_{\text{pa}_D(\tau)}) \log p(i_\tau | i_{\text{pa}_D(\tau)}) \tag{4.3}$$

separately for each $\tau \in \mathcal{T}$. Combining the optima from these separate constrained maximizations according to (3.2) yields a maximum $\hat{p}$ of the likelihood function over $\mathbf{P}_{\text{IV}}(G)$. It is an open question whether or not the likelihood function of the model $\mathbf{P}_{\text{IV}}(G)$ can be multimodal.

In some situations, some of the components of a maximum likelihood estimate in the model $\mathbf{P}_{\text{IV}}(G)$ may be empirical proportions. Recall that a set of vertices $\alpha \subseteq V$ is *complete* if every pair of vertices in $\alpha$ is joined by an edge.

**Proposition 12.** (i) *If the chain component $\tau$ of the chain graph $G$ is a complete set and $\text{pa}_D(\tau) = \varnothing$, then the maximum likelihood estimator of the marginal probability $p(i_\tau)$ is the empirical proportion*

$$\hat{p}(i_\tau) = \frac{n(i_\tau)}{n}.$$

(ii) *If the chain component $\tau$ is a singleton, then the maximum likelihood estimator of the conditional probability $p(i_\tau | i_{\text{pa}_D})$ is the empirical proportion*

$$\hat{p}(i_\tau | i_{\text{pa}_D}) = \frac{n(i_\tau, i_{\text{pa}_D(\tau)})}{n(i_{\text{pa}_D})}.$$

**Proof.** Both observations are immediate consequences of the likelihood factorization provided by (4.2) and (3.8). □

The Möbius parameters of a model $\mathbf{P}_{\text{IV}}(G)$ generally satisfy nonlinear polynomial inequalities. However, if the chain components of $G$ are complete, then linear structure arises, which leads to the following fact:

**Theorem 13.** *If the chain component $\tau$ of the chain graph $G$ is complete and all counts $n(i_\tau, i_{\text{pa}_D(\tau)})$ are positive then the component log-likelihood function $\ell_\tau$ in the model $\mathbf{P}_{\text{IV}}(G)$ (see (4.3)) has a unique local and thus global maximum.*

**Proof.** Under the assumed completeness of $\tau$, Theorem 8(ii) imposes no constraints on the probabilities $p(i_\tau|i_{\text{pa}_D(\tau)})$. Since condition (iii) in Theorem 8 imposes only linear constraints, each $p(i_\tau|i_{\text{pa}_D(\tau)})$ is a linear function of the Möbius parameters $q_\gamma(j_\gamma|i_{\text{pa}_G(\gamma)})$, $\gamma \in \mathcal{C}(\tau)$, defined in (3.7). If $f$ denotes the injective linear map from these Möbius parameters to the probabilities $p(i_\tau|i_{\text{pa}_D(\tau)})$, then the log-likelihood function in terms of the Möbius parameters is strictly concave because it is the composition of $f$ and the strictly concave function $\ell_\tau$ in (4.3). Moreover, the domain of definition of this log-likelihood function is the interior of a polyhedron, and in particular convex. In order to see this note that the probabilities $p(i_\tau|i_{\text{pa}_D(\tau)})$ lie in a Cartesian product of open probability simplices. The preimage of this Cartesian product under the linear map $f$ is the interior of a polyhedron and gives the desired domain of definition; see also [24], Section 1.2.1. The claim now follows because a strictly concave function has a unique local maximum over a convex set. □

We remark that Theorem 13 also applies to discrete chain graph models of type II (AMP) because $\mathbf{P}_{\text{II}}(G) = \mathbf{P}_{\text{IV}}(G)$ if $G$ has complete chain components.

## 5. A non-smooth model of type II (AMP)

Prior work on discrete models of type I (LWF) and our new results on type IV models establish that both these model classes comprise only smooth models. In this section we give an example that shows that the same does not hold for discrete chain graph models of type II (AMP).

Let $G = (V, E)$ be the chain graph in Figure 3. The model $\mathbf{P}_{\text{II}}(G)$ contains the positive distributions for which

$$2 \perp\!\!\!\perp 4 | \{1, 3\} \quad \text{and} \quad 1 \perp\!\!\!\perp \{2, 4\};$$

recall Example 3. In order to analyze $\mathbf{P}_{\text{II}}(G)$ we exploit that $\mathbf{P}_{\text{II}}(G) \subset \mathbf{P}_{\text{IV}}(\bar{G})$, where $\bar{G}$ is the chain graph obtained by adding the edge $2 - 4$ to $G$. The model $\mathbf{P}_{\text{IV}}(\bar{G})$ comprises the positive distributions satisfying $1 \perp\!\!\!\perp \{2, 4\}$ and, according to (3.7), this model can be parametrized using the marginal probabilites

$$q_1(j_1), \qquad q_2(j_2), \qquad q_4(j_4), \qquad q_{24}(j_2, j_4) \tag{5.1}$$

and the conditional probabilities

$$q_3(j_3|i_1), \qquad q_{23}(j_2, j_3|i_1), \qquad q_{34}(j_3, j_4|i_1), \qquad q_{234}(j_2, j_3, j_4|i_1), \tag{5.2}$$

where $i_1 \in [d_1]$ and $j_k \in [d_k - 1]$ for $k = 2, 3, 4$. The conditional independence $2 \perp\!\!\!\perp 4 | \{1, 3\}$, however, is not exhibited by a generic distribution in $\mathbf{P}_{\text{IV}}(\bar{G})$ and leads to additional constraints.

For $i_1 \in [d_1]$ and $i_3 \in [d_3 - 1]$, let $Q^{(i_1,i_3)}$ be the $d_2 \times d_4$-matrix that has entries

$$Q^{(i_1,i_3)}_{i_2 i_4} = \begin{cases} q_{234}(i_2, i_3, i_4|i_1), & \text{if } i_2 < d_2 \text{ and } i_4 < d_4, \\ q_{23}(i_2, i_3|i_1), & \text{if } i_2 < d_2 \text{ and } i_4 = d_4, \\ q_{34}(i_3, i_4|i_1), & \text{if } i_2 = d_2 \text{ and } i_4 < d_4, \\ q_3(i_3|i_1), & \text{if } i_2 = d_2 \text{ and } i_4 = d_4. \end{cases} \quad (5.3)$$

For $i_1 \in [d_1]$ and $i_3 = d_3$, we define $Q^{(i_1,d_3)}$ to be the matrix with entries

$$Q^{(i_1,d_3)}_{i_2 i_4} = \begin{cases} q_{24}(i_2, i_4) - q_{234}(i_2, +, i_4|i_1), & \text{if } i_2 < d_2 \text{ and } i_4 < d_4, \\ q_2(i_2) - q_{23}(i_2, +|i_1), & \text{if } i_2 < d_2 \text{ and } i_4 = d_4, \\ q_4(i_4) - q_{34}(+, i_4|i_1), & \text{if } i_2 = d_2 \text{ and } i_4 < d_4, \\ 1 - q_3(+|i_1), & \text{if } i_2 = d_2 \text{ and } i_4 = d_4. \end{cases} \quad (5.4)$$

In (5.4), the replacement of index $i_3$ by $+$ stands for summation over $i_3 \in [d_3 - 1]$ such that, for example,

$$q_3(+|i_1) = \sum_{i_3=1}^{d_3-1} q_3(i_3|i_1).$$

**Proposition 14.** *Let $G$ be the graph in Figure 3, and $p \in \mathbf{P}_{\text{II}}(\bar{G}) = \mathbf{P}_{\text{IV}}(\bar{G})$. Then, $p \in \mathbf{P}_{\text{II}}(G)$ if and only if the Möbius parameters from (5.1) and (5.2) satisfy that for all $i_1 \in [d_1]$ and $i_3 \in [d_3]$, the matrix $Q^{(i_1,i_3)}$ has a rank of one at most.*

**Proof.** A joint distribution with probability vector $p$ satisfies $2 \perp\!\!\!\perp 4 | \{1, 3\}$ if and only if for all $i_1 \in [d_1]$ and $i_3 \in [d_3]$, the $d_2 \times d_4$-matrix

$$P^{(i_1,i_3)}_{2\perp 4|\{1,3\}} = \begin{pmatrix} p(i_1, 1, i_3, 1) & \cdots & p(i_1, 1, i_3, d_4) \\ \vdots & & \vdots \\ p(i_1, d_2, i_3, 1) & \cdots & p(i_1, d_2, i_3, d_4) \end{pmatrix} \quad (5.5)$$

has a rank of one at most; see, for example, [24], Section 1.5. Using Lemma 7 and Theorem 8, each matrix $P^{(i_1,i_3)}_{2\perp 4|\{1,3\}}$ can be rewritten in terms of the Möbius parameters in (5.1) and (5.2). This requires forming polynomial expressions in the Möbius parameters, but because $\{2, 3, 4\}$ is a complete set in $\bar{G}$, these expressions are equal to the product of a linear term and a marginal probability $q_1(i_1) = P(X_1 = i_1)$. Since $p$ is positive, we can cancel out the marginal probability arriving at a matrix filled with linear expressions; this is equivalent to conditioning on variable $X_1$. After adding rows 1 to $d_2 - 1$ to the last row with index $d_2$ and columns 1 to $d_4 - 1$ to column $d_4$, we arrive at the matrix $Q^{(i_1,i_3)}$. These row and column operations preserve rank and thus the claim follows. □

In order to make our point about non-smoothness of type II (AMP) models, we consider four binary variables, that is, $d_1 = d_2 = d_3 = d_4 = 2$. In this case, there are twelve Möbius parameters in (5.1) and (5.2), and the rank constraints in Proposition 14 require the vanishing of four $2 \times 2$-determinants. For more compact notation, let

$$q_\alpha = q_\alpha(1_\alpha), \quad \alpha \subseteq \{2,4\}$$

and

$$q_{\alpha|i} = q_\alpha(1_\alpha|i), \quad \{3\} \subseteq \alpha \subseteq \{2,3,4\}.$$

Then the two determinants for $i_3 = 1$ yield the equations

$$q_{23|i}q_{34|i} = q_{3|i}q_{234|i}, \qquad i = 1,2, \tag{5.6}$$

whereas the two determinants for $i_3 = d_3 = 2$ yield

$$q_{3|i}q_{24} - q_{23|i}q_4 - q_{34|i}q_2 + q_{234|i} = q_{24} - q_2q_4, \qquad i=1,2. \tag{5.7}$$

We remark that the equations in (5.6) are an instance of a factorization in undirected graphical model [21]: the singleton $\{3\}$ is a separator of the two cliques $\{2,3\}$ and $\{3,4\}$ in the undirected induced subgraph $G_{\{2,3,4\}}$.

The four equations in (5.6) and (5.7) define an eight-dimensional algebraic set in $\mathbb{R}^{11}$; we omit the irrelevant Möbius parameter $q_1$. Using the software Singular [15] we can compute the singularities of this set. (See [3] for a definition of singularities.) We find that the singular locus is determined by the equations

$$q_2q_4 = q_{24}, \qquad q_{3|i}q_2 = q_{23|i}, \qquad q_{3|i}q_4 = q_{34|i}, \qquad q_2q_{3|i}q_4 = q_{234|i}, \qquad i=1,2, \tag{5.8}$$

which by an appeal to Theorem 8 implies the following fact:

**Proposition 15.** *Let $G$ be the graph in Figure 3, and $G_{\text{sing}}$ the subgraph that has the edges $2-3$ and $3-4$ deleted. If all variables are binary ($d_i = 2$), then the singular locus of $\mathbf{P}_{\text{II}}(G)$ is equal to the submodel $\mathbf{P}_{\text{II}}(G_{\text{sing}}) = \mathbf{P}_{\text{IV}}(G_{\text{sing}})$.*

An example of a statistical consequence of the non-smoothness of $\mathbf{P}_{\text{II}}(G)$ is that a $\chi^2$-approximation is inappropriate for the likelihood ratio test of $\mathbf{P}_{\text{II}}(G_{\text{sing}})$ versus $\mathbf{P}_{\text{II}}(G)$; compare [6].

*Remark 16.* We conjecture that $\mathbf{P}_{\text{II}}(G)$ is non-smooth regardless of the number of levels $d_i$ for the four random variables. Using Singular, we were able to verify the claim of Proposition 15 when $X_3$ is ternary, that is, $d_1 = d_2 = d_4 = 2$ and $d_3 = 3$. Moreover, we could compute the case $d_2 = d_3 = d_4 = 2$ and $d_1 = 3$ for which $\mathbf{P}_{\text{II}}(G_{\text{sing}})$ is a only a proper subset of the singular locus.

## 6. A non-smooth model of type III

Let $G$ be the chain graph in Figure 3. As seen in the previous section, the binary type II model $\mathbf{P}_{\mathrm{II}}(G)$ is non-smooth. Nevertheless one can give a rational parametrization of $\mathbf{P}_{\mathrm{II}}(G)$ by solving the equations (5.6) and (5.7). In this section we show that the binary type III model $\mathbf{P}_{\mathrm{III}}(G)$ is the union of two strict submodels defined by polynomial equations, which implies that the model is non-smooth and cannot be parametrized using a rational map. The non-existence of a rational parametrization follows from the fact that an algebraic set with rational parametrization is irreducible [3], Section 4.5. An algebraic set, that is, a set defined by polynomial equations, is irreducible if it cannot be decomposed into a finite union of strict algebraic subsets.

Let $X_1, \ldots, X_4$ be binary variables in correspondence with the nodes of $G$. As stated in Example 3, $\mathbf{P}_{\mathrm{III}}(G)$ comprises the positive distributions for which

$$2 \perp\!\!\!\perp 4 | 1 \quad \text{and} \quad 1 \perp\!\!\!\perp \{2,4\} | 3. \tag{6.1}$$

Define the two $2 \times 2$-matrices

$$P_{2 \perp\!\!\!\perp 4 | 1}^{(i)} = \begin{pmatrix} p_{i1+1} & p_{i1+2} \\ p_{i2+1} & p_{i2+2} \end{pmatrix}, \qquad i = 1, 2, \tag{6.2}$$

and the two $2 \times 4$-matrices

$$P_{1 \perp\!\!\!\perp \{2,4\} | 3}^{(k)} = \begin{pmatrix} p_{11k1} & p_{11k2} & p_{12k1} & p_{12k2} \\ p_{21k1} & p_{21k2} & p_{22k1} & p_{22k2} \end{pmatrix}, \qquad k = 1, 2. \tag{6.3}$$

Here, $p_{ijk\ell} = p(i,j,k,\ell)$ and $p_{ij+\ell} = P(X_1 = i, X_2 = j, X_4 = \ell)$ for $i, j, \ell = 1, 2$.

A probability distribution on $[2]^4$ satisfies the two conditional independences in (6.1) if and only if the four matrices in (6.2) and (6.3) have a rank of one at most. This rank condition together with the constraint that the probabilities $p_{ijk\ell}$ sum to one defines a seven-dimensional algebraic set in $\mathbb{R}^{16}$. By computing a primary decomposition using Singular, this set is seen to decompose into the union of two strict algebraic subsets that both have dimension seven. (See again [3] for an introduction to the involved algebraic concepts.)

**Proposition 17.** *Let $G$ be the graph in Figure 3. If all variables are binary $(d_i = 2)$, then a positive probability vector $p = (p_{ijk\ell})$ is in $\mathbf{P}_{\mathrm{III}}(G)$ if and only if at least one of the following two conditions is met:*

(i) $1 \perp\!\!\!\perp \{2,3,4\}$ *and* $2 \perp\!\!\!\perp \{1,4\}$, *or*
(ii) $2 \perp\!\!\!\perp 4 | 1$ *and* $1 \perp\!\!\!\perp \{2,4\} | 3$ *and*

$$p_{1121} p_{2222} - p_{1122} p_{2221} = p_{1111} p_{2212} - p_{1112} p_{2211}, \tag{6.4}$$

$$p_{1221} p_{2122} - p_{1222} p_{2121} = p_{1211} p_{2112} - p_{1212} p_{2111}. \tag{6.5}$$

The equations (6.4) and (6.5) are equalities of $2 \times 2$-minors of the two matrices

$$P^{(k)}_{\{1,2\} \perp\!\!\!\perp 4|3} = \begin{pmatrix} p_{11k1} & p_{12k1} & p_{21k1} & p_{22k1} \\ p_{11k2} & p_{12k2} & p_{21k2} & p_{22k2} \end{pmatrix}, \qquad k=1,2. \tag{6.6}$$

associated with the conditional independence $\{1,2\} \perp\!\!\!\perp 4|3$. Hence, for these equations to hold some homogeneity between possible conditional dependence between $(X_1, X_2)$ and $X_4$ given $X_3 = 1$ and possible conditional dependence between $(X_1, X_2)$ and $X_4$ given $X_3 = 2$ is required.

In order to show that the two components of $\mathbf{P}_{\text{III}}(G)$ are indeed non-trivial and distinct from each other over the interior of the probability simplex we give the following two examples: The probability vector

$$\begin{pmatrix} p_{1111} & p_{1112} & p_{1121} & p_{1122} \\ p_{1211} & p_{1212} & p_{1221} & p_{1222} \\ p_{2111} & p_{2112} & p_{2121} & p_{2122} \\ p_{2211} & p_{2212} & p_{2221} & p_{2222} \end{pmatrix} = \frac{1}{32} \begin{pmatrix} 1 & 3 & 3 & 1 \\ 3 & 1 & 1 & 3 \\ 1 & 3 & 3 & 1 \\ 3 & 1 & 1 & 3 \end{pmatrix} \tag{6.7}$$

satisfies the condition in Proposition 17(i) but not the one in (ii), whereas

$$\begin{pmatrix} p_{1111} & p_{1112} & p_{1121} & p_{1122} \\ p_{1211} & p_{1212} & p_{1221} & p_{1222} \\ p_{2111} & p_{2112} & p_{2121} & p_{2122} \\ p_{2211} & p_{2212} & p_{2221} & p_{2222} \end{pmatrix} = \frac{1}{32} \begin{pmatrix} 2 & 2 & 2 & 2 \\ 2 & 2 & 1 & 1 \\ 3 & 3 & 2 & 2 \\ 3 & 3 & 1 & 1 \end{pmatrix} \tag{6.8}$$

satisfies (ii) but not (i).

**Remark 18.** More complicated computations in `Singular` are feasible. If $X_3$ is ternary ($d_1 = d_2 = d_4 = 2$ and $d_3 = 3$), then the algebraic set corresponding to the model $\mathbf{P}_{\text{III}}(G)$ breaks into six components. If $X_1$ is the only ternary variable ($d_2 = d_3 = d_4 = 2$ and $d_1 = 3$), then the set is irreducible. Nevertheless, we conjecture that $\mathbf{P}_{\text{III}}(G)$ is non-smooth regardless of the choice of $d_1, \ldots, d_4$.

## 7. Discussion

The main contribution of this paper concerns discrete chain graph models of type IV. These models are related to multivariate regression and can also be derived from a path diagram interpretation for the chain graph. In fact, the block-recursive Markov property of type IV can be shown to be equivalent to the global Markov property discussed, for instance, in [20, 25, 26]. In this paper we showed that, just like their Gaussian analogs, type IV models are curved exponential families (Corollary 10). This brings with it all the convenience of the standard asymptotic theory for likelihood inference.

Practical use of the discrete models requires algorithms for maximization of the likelihood function, and at least two approaches are possible. On one hand, one may write the likelihood function as a function of the Möbius parameters from (3.7) and then apply general optimizers. On the other hand, the iterative conditional fitting algorithm of [8] can be extended to the chain graph case, as explained in [5].

The approach we took when studying type IV models is also useful for analyzing discrete models of type II (AMP). While much of the related work on models for categorical data employs log-linear expansions (see e.g. [2, 19, 23]), the change of coordinates in Theorem 8 remains at the level of conditional probabilities. This preserves the algebraic structure that was useful for showing that the type II class includes models with singularities (Section 5). An interesting problem for future research would be to characterize all chain graphs that yield smooth discrete (or perhaps more concretely, binary) chain graph models of type II. In particular, an interesting question is whether there exist chain graphs $G$ for which $\mathbf{P}_{II}(G)$ is smooth and such that there does not exist a chain graph $\bar{G}$ with $\mathbf{P}_{II}(G) = \mathbf{P}_{I}(\bar{G})$ or $\mathbf{P}_{II}(G) = \mathbf{P}_{IV}(\bar{G})$. Similar questions arise for models of type III that may also be non-smooth (Section 6).

Finally, we recall Remark 5, where we commented on possible pairwise Markov interpretations of chain graphs. While these can be awkward in the sense that the focus may be on distributions for which $v \perp\!\!\!\perp u$, $v \perp\!\!\!\perp w$ but not $v \perp\!\!\!\perp \{u, w\}$, it is interesting that the pairwise type II interpretation of the graph in Figure 3 yields a smooth model.

# Acknowledgement

Research partially supported by National Science Foundation Grants DMS-0505612 and 0746265.